\documentclass[11pt,letterpaper,reqno]{amsart}

\usepackage{times}
\usepackage[T1]{fontenc}
\usepackage{mathrsfs}
\usepackage{latexsym}
\usepackage[dvips]{graphics}
\usepackage{epsfig}
\usepackage{amsmath,amsfonts,amsthm,amssymb,amscd}
\input amssym.def
\input amssym.tex

\newcommand\be{\begin{equation}}
\newcommand\ee{\end{equation}}
\newcommand\bea{\begin{eqnarray}}
\newcommand\eea{\end{eqnarray}}
\newcommand\bi{\begin{itemize}}
\newcommand\ei{\end{itemize}}
\newcommand\ben{\begin{enumerate}}
\newcommand\een{\end{enumerate}}
\newcommand\bc{\begin{center}}
\newcommand\ec{\end{center}}
\newcommand\ba{\begin{array}}
\newcommand\ea{\end{array}}
\newcommand{\foh}{\frac{1}{2}}  %onehalf
\newtheorem{thm}{Theorem}[section]

\newtheorem{lem}[thm]{Lemma}

\theoremstyle{definition}

\begin{document}

\title{A Probabilistic Proof of Wallis's Formula for $\pi$}

\author{Steven J. Miller}


\date{\today}

%\begin{abstract}
%\end{abstract}

%\pagebreak

\maketitle

%\setcounter{equation}{0}
    %=======================================================
    %   You need the setcounter command above to reset the
    %   equation counting in each section. without it, you'd
    %   have 1.1, 1.2, 1.3, 2.4, 2.5 -- with it (you put it
    %   after each new section) you have 1.1, 1.2, 1.3, 2.1,
    %========================================================

%%%%%%%%%%%%%%%%%%%%%%%%%%%%%%%%%%%%%%%%%%%%%%%%%%%%%%%%%%%%%%%%%%%%%%%%%%%%%%%%%%%%%%%%%%%%%
%%%%%%%%%%%%%%%%%%%%%%%%%%%%%%%%%%%%%%%%%%%%%%%%%%%%%%%%%%%%%%%%%%%%%%%%%%%%%%%%%%%%%%%%%%%%%
%\setcounter{equation}{0}

There are many beautiful formulas for $\pi$ (see for example
\cite{BB}). The purpose of this note is to introduce an alternate
derivation of Wallis's product formula, equation \eqref{eq:wallis},
which could be covered in a first course on probability, statistics,
or number theory. We quickly review other famous formulas for $\pi$,
recall some needed facts from probability, and then derive Wallis's
formula. We conclude by combining some of the other famous formulas
with Wallis's formula to derive an interesting expression for $\log
(\pi/2)$ (equation \eqref{eq:logpiformula}).

Often in a first-year calculus course students encounter the
Gregory-Leibniz formula, \be \frac{\pi}4 \ = \ 1 - \frac13 + \frac15
- \frac17 + \cdots \ = \ \sum_{n=0}^\infty \frac{(-1)^n}{2n+1}.
\nonumber\ \ee The proof uses the fact that the derivative of
$\arctan x$ is $1/(1+x^2$), so $\pi/4 = \int_0^1 dx/(1+x^2)$. To
complete the proof, expand the integrand with the geometric series
formula and then justify interchanging the order of integration and
summation.

Another interesting formula involves Bernoulli numbers and the
Riemann zeta function. The Bernoulli numbers $B_k$ are the
coefficients in the Taylor series \be \frac{t}{e^t - 1} \ = \ 1 -
\frac{t}2 + \sum_{k=2}^\infty \frac{B_k t^k}{k!}; \nonumber \ee each
$B_k$ is rational. The Riemann zeta function is $\zeta(s) =
\sum_{n=1}^\infty n^{-s}$, which converges for real part of $s$
greater than $1$. Using complex analysis one finds (see for instance
[10, p. 365] or [18, pp. 179--180]) that \be\label{eq:zetabernpi}
\zeta(2k) \ = \ - \frac{(-4)^k B_{2k}}{2\cdot 2k!} \cdot \pi^{2k},
\nonumber\ \ee yielding formulas for $\pi$ to any even
power.\footnote{An amusing consequence of these formulas is a proof
of the infinitude of primes. Using unique factorization, one can
show that $\zeta(s)$ also equals $\prod_p (1-p^{-s})^{-1}$, where
$p$ runs over all primes. As $\pi^2$ is irrational and $\zeta(2) =
\pi^2/6$, there must be infinitely many primes: if there were only
finitely many then $\pi^2/6 = \prod_p (1-p^{-2})^{-1}$ would be
rational! See \cite{MSW} for explicit lower bounds on $\pi(x)$
derivable from upper bounds for the irrationality measure of
$\zeta(2)$, and \cite{MT-B} for more details on the numerous
connections between $\zeta(s)$ and number theory.} In particular,
$\pi^2/6 = \sum_n n^{-2}$ and $\pi^4/90 = \sum_n n^{-4}$.

One of the most interesting formulas for $\pi$ is a multiplicative
one due to Wallis (1665): \be\label{eq:wallis} \frac{\pi}2\ \ = \ \
\frac{2\cdot 2}{1\cdot 3}\ \ \frac{4\cdot 4}{3\cdot 5}\ \
\frac{6\cdot 6}{5\cdot 7}\ \ \frac{8\cdot 8}{7\cdot 9}\ \cdots\ \ =
\ \ \prod_{n=1}^\infty \frac{2n \cdot 2n}{(2n-1)(2n+1)}.  \ee Common
proofs use the infinite product expansion for $\sin x$ (see [18, p.
142]) or induction to prove formulas for integrals of powers of
$\sin x$ (see [3, p. 115]). We present a mostly elementary proof
using standard facts about probability distributions encountered in
a first course on probability or statistics (and hence the
title).\footnote{For a statistical proof involving an experiment and
data, see the chapter on Buffon's needle in \cite{AZ} (page 133): if
you have infinitely many parallel lines $d$ units apart, then the
probability that a ``randomly'' dropped rod of length $\ell \le d$
crosses one of the lines is $2\ell/\pi d$. Thus you can calculate
$\pi$ by throwing many rods on the grid and counting the number of
intersections.} The reason we must write ``mostly elementary'' is
that at one point we appeal to the Dominated Convergence Theorem. It
is possible to bypass this and argue directly, and we sketch the
main ideas for the interested reader.

Recall that a continuous function $f(x)$ is a continuous probability
distribution if (1) $f(x) \ge 0$ and (2) $\int_{-\infty}^\infty f(x)
dx = 1$. We immediately see that if $g(x)$ is a nonnegative
continuous function whose integral is finite then there exists an $a
> 0$ such that $a g(x)$ is a continuous probability distribution
(take $a = 1 / \int_{-\infty}^\infty g(x)dx$). This simple
observation is a key ingredient in our proof, and is an extremely
important technique in mathematics; the proof of Wallis's formula is
but one of many applications.\footnote{A nice application of
Wallis's formula is in determining the universal constant in
Stirling's formula for $n!$; see \cite{Sc} for some history and
applications.} In fact, this observation greatly simplifies numerous
calculations in random matrix theory, which has successfully modeled
diverse systems ranging from energy levels of heavy nuclei to the
prime numbers; see \cite{Con,MT-B} for introductions to random
matrix theory and \cite{Meh} for applications of this technique to
the subject. One of the purposes of this paper is to introduce
students to the consequences of this simple observation.

Our proof relies on two standard functions from probability, the
Gamma function and the Student $t$-distribution. The Gamma function
$\Gamma(x)$ is defined by \be \Gamma(x) \ = \ \int_0^\infty e^{-t}
t^{x-1}dt. \nonumber\ \ee Note that this integral is well defined if
the real part of $x$ is positive. Integrating by parts yields
$\Gamma(x+1) = x\Gamma(x)$. This implies that if $n$ is a
nonnegative integer then $\Gamma(n+1) = n!$; thus the Gamma function
generalizes the factorial function (see \cite{Sr} for more on the
Gamma function, including another proof of Wallis's formula
involving the Gamma function). We need the following:

\noindent\textbf{Claim: $\Gamma(1/2) = \sqrt{\pi}$.}

\begin{proof} In the integral for $\Gamma(1/2)$, change variables by
setting $u = \sqrt{t}$ (so $dt=2udu = 2\sqrt{t}du$). This yields \be
\Gamma(1/2) \ = \ 2 \int_0^\infty e^{-u^2} du \ = \
\int_{-\infty}^\infty e^{-u^2}du. \nonumber\ \ee This integral is
well-known to equal $\sqrt{\pi}$ (see page 542 of \cite{Ar}). The
standard proof is to square the integral and convert to polar
coordinates: \bea \Gamma(1/2)^2 & \ = \ &  \int_{-\infty}^\infty
e^{-u^2} du \int_{-\infty}^\infty e^{-v^2} dv \ = \  \int_0^\infty
\int_0^{2\pi} e^{-r^2} r dr d\theta \ = \ \pi. \nonumber\ \eea
\end{proof}

In fact, our proof above shows \be\label{eq:defstdnormal}
\int_{-\infty}^\infty \frac{1}{\sqrt{2\pi}}\ e^{-t^2/2} dt \ = \ 1.
\ee This density is called the standard normal (or Gaussian). This
is one of the most important probability distributions, and we shall
see it again when we look at the Student $t$-distribution. If $g$ is
a continuous probability density, then we say that the random
variable $Y$ has distribution $g$ if for any interval $[a,b]$ the
probability that $Y$ takes on a value in $[a,b]$ is $\int_a^b
g(y)dy$. The celebrated Central Limit Theorem (see [6, p. 515] for a
proof) states that for many continuous densities $g$, if $X_1,
\dots, X_n$ are independent random variables, each with density $g$,
then as $n\to \infty$ the distribution of $(Y_n - \mu)/(\sigma/\sqrt{n})$
converges to the standard normal (where $Y_n = (X_1 + \cdots +
X_n)/n$ is the sample average, $\mu$ is the mean of $g$, and
$\sigma$ is its standard deviation\footnote{The mean $\mu$ of a
distribution is its average value: $\mu = \int x g(x)dx$. The
standard deviation $\sigma$ measures how spread out a distribution
is about its average value: $\sigma^2 = \int (x-\mu)^2 g(x)dx$.}).

The second function we need is the Student\footnote{Student was the
pen name of William Gosset.} $t$-distribution (with $\nu$ degrees of
freedom): \be f_\nu(t) \ = \
\frac{\Gamma\left(\frac{\nu+1}2\right)}{\sqrt{\pi\nu}\
\Gamma\left(\frac{\nu}2\right)} \ \cdot \ \left(1 +
\frac{t^2}{\nu}\right)^{-\frac{\nu+1}2} \ = \ c_\nu \left(1 +
\frac{t^2}{\nu}\right)^{-\frac{\nu+1}2}; \nonumber\ \ee here
$\nu$ is a positive integer and $t$ is any real number.\\

\noindent\textbf{Claim: The Student $t$-distribution is a continuous
probability density.}

\begin{proof} As $f_\nu(t)$ is clearly continuous and nonnegative, to show $f_\nu(t)$ is a probability
density it suffices to show that it integrates to $1$. We must
therefore show that \be \int_{-\infty}^\infty \left(1 +
\frac{t^2}{\nu}\right)^{-\frac{\nu+1}2}dt \ = \ \frac{\sqrt{\pi\nu}\
\Gamma\left(\frac{\nu}2\right)}{\Gamma\left(\frac{\nu+1}2\right)}.
\nonumber\ \ee As the integrand is symmetric, we may integrate from
$0$ to infinity and double the result. Letting $t = \sqrt{\nu} \tan
\theta$ (so $dt = \sqrt{\nu} \sec^2\theta d\theta$) we find \be
\int_{-\infty}^\infty \left(1 +
\frac{t^2}{\nu}\right)^{-\frac{\nu+1}2}dt \ = \ 2\sqrt{\nu}
\int_0^{\pi/2} \frac{\sec^2 \theta d\theta}{\sec^{\nu+1}\theta} \ =
\ 2\sqrt{\nu} \int_0^{\pi/2} \cos^{\nu-1} \theta d\theta. \nonumber\
\ee The proof follows immediately from two properties of the Beta
function (see [2, p. 560]): \bea B(p,q) & \ =\ & \Gamma(p)\Gamma(q)
/ \Gamma(p+q) \nonumber\\ B(m+1,n+1) &=& 2\int_0^{\pi/2} \cos^{2m+1}
(\theta) \sin^{2n+1}(\theta) d\theta; \eea an elementary proof
without appealing to properties of the Beta function is given in
Appendix \ref{sec:appendix}.
\end{proof}

The Student $t$-distribution arises in statistical analyses where
the sample size $\nu$ is small and each observation is normally
distributed with the same mean and the same (unknown) variance (see
\cite{Go,MM}). The reason the Student $t$-distribution is used only
for small samples sizes is that as $\nu \to \infty$, $f_\nu(t)$
converges to the standard normal; proving this will yield Wallis's
formula. While we can prove this by invoking the Central Limit
Theorem, we may also see this directly by recalling that \be e^x \ =
\ \lim_{N\to\infty} \left(1 + \frac{x}{N}\right)^N. \nonumber\ \ee
We therefore have \be \lim_{\nu \to\infty}\left(1 +
\frac{t^2}{\nu}\right)^{-\nu/2} \ = \ \left(e^{t^2}\right)^{-1/2} \
= \ e^{-t^2/2}. \nonumber\ \ee As $f_\nu(t)$ is a probability
distribution for all positive integers $\nu$, it integrates to $1$
for all such $\nu$, which is equivalent to \be \frac1{c_\nu} \ = \
\frac{\sqrt{\pi\nu}\
\Gamma\left(\frac{\nu}2\right)}{\Gamma\left(\frac{\nu+1}2\right)} \
= \ \int_{-\infty}^\infty \left(1 +
\frac{t^2}{\nu}\right)^{-\frac{\nu+1}2}dt.\nonumber\ \ee Taking the
limit as $\nu\to\infty$ yields \bea \lim_{\nu \to \infty}
\frac1{c_\nu} & \ = \ & \lim_{\nu \to \infty} \int_{-\infty}^\infty
\left(1 + \frac{t^2}{\nu}\right)^{-\frac{\nu+1}2}dt \nonumber\\ & =
& \int_{-\infty}^\infty \lim_{\nu \to \infty} \left(1 +
\frac{t^2}{\nu}\right)^{-\frac{\nu+1}2}dt \ = \
\int_{-\infty}^\infty e^{-t^2/2} dt \ = \ \sqrt{2\pi}. \nonumber\
\eea Some work is necessary of course to justify interchanging the
integral and the limit; this justification is why our argument is
only ``mostly elementary''. A standard proof uses the Dominated
Convergence Theorem (see [7, p. 54] or [9, p. 238]) to show that as
$\nu\to\infty$ the $t$-distribution converges to the standard
normal;\footnote{For completeness we sketch how such an argument
could proceed. If $t \in [-\log^2 \nu, \log^2 \nu]$ then
$|(1+t^2/\nu)^{-(\nu+1)/2} - \exp(-t^2/2)|$  tends to zero rapidly
with $\nu$. Further, if $f(t)$ is the density of the standard
normal, then $\int_{|t| \ge \log^2 \nu} f(t)dt$ and $\int_{|t| \ge
\log^2 \nu} f_\nu(t)dt$ also tend to zero rapidly with $\nu$.
Careful bookkeeping shows that the normalization constants $c_\nu$
must therefore approach $c=1/\sqrt{2\pi}$.} one may take $2008
\exp(-t^2/2008)$ as the dominating function. We have therefore shown
that \be\label{eq:climitcnugamma} c\ = \ \lim_{\nu \to \infty} c_\nu
\ = \ \lim_{\nu \to \infty}
\frac{\Gamma\left(\frac{\nu+1}2\right)}{\sqrt{\pi\nu}\
\Gamma\left(\frac{\nu}2\right)} \ = \ \frac1{\sqrt{2\pi}}. \ee The
fact that $c=1/\sqrt{2\pi}$ is the key step in our proof of Wallis's
formula. We have calculated the limit by using our observation that
a probability distribution must integrate to $1$; calculating it by
brute force
analysis of the Gamma factors yields our main result.\\

\noindent\textbf{Theorem: Wallis's formula is true.}

\begin{proof} The proof follows from expanding the Gamma functions and substituting
into \eqref{eq:climitcnugamma}; we highlight the main steps. Let
$\nu = 2m$. Using $\Gamma(n+1) = n\Gamma(n)$ and $\Gamma(1/2) =
\sqrt{\pi}$ we find that \bea \Gamma\left(\frac{2m+1}2\right) & \ =
\ & (2m-1) (2m-3) \cdots 5 \cdot 3 \cdot 1 \cdot \sqrt{\pi} / 2^m.
\nonumber\ \eea As $\Gamma(m) = (m-1)!$, after some algebra we find
that \bea c_{2m} \ = \ \frac{1\cdot 3 \cdot 5 \cdots (2m-3) \cdot
(2m-1)}{2\cdot 4 \cdot 6\ \  \cdots\  \ (2m-2)\ \ \ \cdot 2m \ \ \ \
} \ \frac{\sqrt{m}}{\sqrt{2}}. \nonumber\ \eea Multiplying by
$1\cdot(2m+1)/(2m+1)$ and regrouping, we find \bea c_{2m} & \ = \ &
\left(\frac{1\cdot 3}{2 \cdot 2} \ \cdot \ \frac{3\cdot 5}{4\cdot
4}\ \cdots\ \frac{(2m-1)(2m+1)}{2m\cdot 2m} \
\frac1{2m+1}\right)^\foh \ \frac{\sqrt{m}}{\sqrt{2}}, \nonumber\
\eea which we rewrite as \bea \prod_{n=1}^m \frac{2n\cdot
2n}{(2n-1)(2n+1)} \ = \ \frac{2\cdot 2}{1\cdot 3}\ \frac{4\cdot
4}{3\cdot 5}\ \cdots\ \frac{2m \cdot 2m}{(2m-1)(2m+1)} \ = \
\frac{m}{(4m+2)c_{2m}^2}. \nonumber\ \eea As $\lim_{m\to\infty}
1/c_{2m}^2 = 2\pi$ and $\lim_{m\to\infty} m/(4m+2) = 1/4$, we have
\be \prod_{n=1}^\infty \frac{2n\cdot 2n}{(2n-1)(2n+1)} \ = \
\frac{\pi}2, \nonumber\ \ee which completes the proof.
\end{proof}

By combining the expansion for $\pi$ from Wallis's formula with
those involving the Bernoulli numbers and the zeta function, we
obtain a proof of the following amusing formula for
$\log (\pi/2)$.\\

\noindent\textbf{Theorem: We have \be\label{eq:logpiformula}
\log\frac{\pi}2\ =\ \sum_{k=1}^\infty \frac{\zeta(2k)}{4^k \cdot k}
\ = \ -\sum_{k=1}^\infty \frac{(-1)^k B_{2k}}{2k\cdot 2k!} \cdot
\pi^{2k}.  \ee}

\begin{proof} The Taylor series of $\log(1-x)$ is
$-\sum_{k=1}^\infty x^k/k$. The $n$th factor in Wallis's formula may
be written as $\left(1 - \frac1{4n^2}\right)^{-1}$. Thus taking
logarithms of Wallis's formula and Taylor expanding yields \be
\log\frac{\pi}2 \ = \ - \sum_{n=1}^\infty \log\left(1 -
\frac1{4n^2}\right) \ = \  \sum_{n=1}^\infty \sum_{k=1}^\infty
\frac{1}{k \cdot (4n^2)^k} \ = \ \sum_{k=1}^\infty \frac{1}{4^k
\cdot k} \sum_{n=1}^\infty\frac{1}{n^{2k}}. \nonumber\ \ee The
$n$-sum gives $\zeta(2k)$, and the claim now follows.
\end{proof}

Note that the above formula for $\log(\pi/2)$ converges well. It is
easy to see that $|\zeta(2k)| \le 2$ (and $\lim_{k\to\infty}
\zeta(2k) = 1$). Thus each additional summand yields at least one
new digit (base $4$). See \cite{So} for additional formulas for
$\log(\pi/2)$.

\ \\

\noindent \textbf{ACKNOWLEDGEMENTS.} The author would like to thank
his statistics classes and participants at $\pi$ Day at Brown for
many enlightening conversations, and Mike Rosen, Ed Scheinerman,
Gopal Srinivasan, and the referees for several useful suggestions
and references. This work was partly supported by NSF
grant DMS0600848.\\

\appendix

%%%%%%%%%%%%%%%%%%%%%%%%%%%%%%%%%%%%%%%%%%%%%%%%%%%%%%%%%%%%%%%%%%%%%%%%
%%%%%%%%%%%%%%%%%%%%%%%%%%%%%%%%%%%%%%%%%%%%%%%%%%%%%%%%%%%%%%%%%%%%%%%%

\section{Elementary calculation of constants}\label{sec:appendix}

For completeness, we provide a more elementary derivation that the
stated constant is the correct normalization constant for the
Student $t$-distribution.

%\begin{lem} $\Gamma(1/2) = \sqrt{\pi}$. \end{lem}
%
%\begin{proof} In the integral for $\Gamma(1/2)$, change variables by
%setting $u = \sqrt{t}$ (so $dt=2udu = 2\sqrt{t}du$). This yields \be
%\Gamma(1/2) \ = \ 2 \int_0^\infty e^{-u^2} du \ = \
%\int_{-\infty}^\infty e^{-u^2}du. \nonumber\ \ee This integral is
%well-known to equal $\sqrt{\pi}$. The standard proof is to square
%the integral and convert to polar coordinates: \bea \Gamma(1/2)^2 &
%\ = \ &  \int_{-\infty}^\infty e^{-u^2} du \int_{-\infty}^\infty
%e^{-v^2} dv \ = \  \int_0^\infty \int_0^{2\pi} e^{-r^2} r dr d\theta
%\ = \ \pi. \nonumber\ \eea
%\end{proof}

\begin{lem} We have \be 2\sqrt{\nu} \int_0^{\pi/2} \cos^{\nu-1} \theta d\theta
\ = \ \frac{\sqrt{\pi\nu}\
\Gamma\left(\frac{\nu}2\right)}{\Gamma\left(\frac{\nu+1}2\right)}.\ee
\end{lem}

\begin{proof} The claim follows by induction; we sketch the main idea.
Assume we have proven the claim for all $\nu \le n$. Then \bea
\int_0^{\pi/2} \cos^{n+1}\theta d\theta & \ = \ & \int_0^{\pi/2} (1
- \sin^2\theta) \cos^{n-1}\theta d\theta \nonumber\\ & = &
\int_0^{\pi/2} \cos^{n-1} \theta d\theta - \int_0^{\pi/2} \sin
\theta \cdot \left(\cos^{n-1} \theta \sin \theta\right) d\theta.
\nonumber\ \eea We integrate the second term on the right by parts,
with $u = \sin \theta$ and $dv = \cos^{n-1} \theta \sin \theta
d\theta$. The $uv$ term vanishes at $0$ and $\pi/2$ and we are left
with \be \int_0^{\pi/2} \cos^{n+1}\theta d\theta \ = \
\int_0^{\pi/2} \cos^{n-1} \theta d\theta - \frac1n \int_0^{\pi/2}
\cos^n\theta \cos\theta d\theta, \nonumber\ \ee which simplifies to
\be \int_0^{\pi/2} \cos^{n+1}\theta d\theta \ = \ \frac{n}{n+1}
\int_0^{\pi/2} \cos^{n-1}\theta d\theta. \nonumber\ \ee The claim
now follows from standard properties of the Gamma function
($\Gamma(m+1) = m\Gamma(m)$ and $\Gamma(1/2) = \sqrt{\pi}$).
\end{proof}

%%%%%%%%%%%%%%%%%%%%%%%%%%%%%%%%%%%%%%%%%%%%%%%%%%%%%%%%%%%%%%%%%%%%%%%%%%%%%%%%%%%%%%%%%%%%%
%%%%%%%%%%%%%%%%%%%%%%%%%%%%%%%%%%%%%%%%%%%%%%%%%%%%%%%%%%%%%%%%%%%%%%%%%%%%%%%%%%%%%%%%%%%%%

\ \\

\noindent Department of Mathematics, Brown University, Providence,
RI 02912 and Department of Mathematics and Statistics, Williams College, Williamstown, MA 01267

\noindent Steven.J.Miller@williams.edu

\end{document}